\documentclass[12pt]{article} 
\usepackage{latexsym} 
\usepackage{amsmath}
\usepackage{mathptmx}
\usepackage{amssymb} 
\usepackage{amssymb} 
\usepackage{amsfonts}
\usepackage{amscd} 
\usepackage[all]{xy}
\usepackage{graphicx}
\usepackage{layout}
\begin{document} 
\newtheorem{Th}{Theorem}[section]
\newtheorem{Cor}{Corollary}[section]
\newtheorem{Prop}{Proposition}[section]
\newtheorem{Lem}{Lemma}[section]
\newtheorem{Def}{Definition}[section]
\newtheorem{Rem}{Remark}[section]
\newtheorem{Ex}{Example}[section]
\newtheorem{stw}{Proposition}[section]


\newcommand{\bet}{\begin{Th}}
\newcommand{\ent}{\stepcounter{Cor}
   \stepcounter{Prop}\stepcounter{Lem}\stepcounter{Def}
   \stepcounter{Rem}\stepcounter{Ex}\end{Th}}


\newcommand{\bec}{\begin{Cor}}
\newcommand{\enc}{\stepcounter{Th}
   \stepcounter{Prop}\stepcounter{Lem}\stepcounter{Def}
   \stepcounter{Rem}\stepcounter{Ex}\end{Cor}}
\newcommand{\bep}{\begin{Prop}}
\newcommand{\enp}{\stepcounter{Th}
   \stepcounter{Cor}\stepcounter{Lem}\stepcounter{Def}
   \stepcounter{Rem}\stepcounter{Ex}\end{Prop}}
\newcommand{\bel}{\begin{Lem}}
\newcommand{\enl}{\stepcounter{Th}
   \stepcounter{Cor}\stepcounter{Prop}\stepcounter{Def}
   \stepcounter{Rem}\stepcounter{Ex}\end{Lem}}
\newcommand{\bef}{\begin{Def}}
\newcommand{\enf}{\stepcounter{Th}
   \stepcounter{Cor}\stepcounter{Prop}\stepcounter{Lem}
   \stepcounter{Rem}\stepcounter{Ex}\end{Def}}
\newcommand{\ber}{\begin{Rem}}
\newcommand{\enr}{
   \stepcounter{Th}\stepcounter{Cor}\stepcounter{Prop}
   \stepcounter{Lem}\stepcounter{Def}\stepcounter{Ex}\end{Rem}}
\newcommand{\bee}{\begin{Ex}}
\newcommand{\ene}{
   \stepcounter{Th}\stepcounter{Cor}\stepcounter{Prop}
   \stepcounter{Lem}\stepcounter{Def}\stepcounter{Rem}\end{Ex}}
\newcommand{\Proof}{\noindent{\it Proof\,}:\ }
\newcommand{\beP}{\Proof}
\newcommand{\enP}{\hfill $\Box$ \par\vspace{5truemm}}

\newcommand{\EE}{\mathbb{E}}
\newcommand{\QQ}{\mathbb{Q}}
\newcommand{\R}{\mathbb{R}}
\newcommand{\C}{\mathbb{C}}
\newcommand{\ZZ}{\mathbb{Z}}
\newcommand{\KK}{\mathbb{K}}
\newcommand{\NN}{\mathbb{N}}
\newcommand{\PP}{\mathbb{P}}
\newcommand{\HH}{\mathbb{H}}
\newcommand{\uuu}{\boldsymbol{u}}
\newcommand{\xxx}{\boldsymbol{x}}
\newcommand{\aaa}{\boldsymbol{a}}
\newcommand{\bbb}{\boldsymbol{b}}
\newcommand{\AAA}{\mathbf{A}}
\newcommand{\BBB}{\mathbf{B}}
\newcommand{\LLL}{\mathbf{L}}
\newcommand{\ccc}{\boldsymbol{c}}
\newcommand{\iii}{\boldsymbol{i}}
\newcommand{\jjj}{\boldsymbol{j}}
\newcommand{\kkk}{\boldsymbol{k}}
\newcommand{\rrr}{\boldsymbol{r}}
\newcommand{\FFF}{\boldsymbol{F}}
\newcommand{\yyy}{\boldsymbol{y}}
\newcommand{\ppp}{\boldsymbol{p}}
\newcommand{\qqq}{\boldsymbol{q}}
\newcommand{\nnn}{\boldsymbol{n}}
\newcommand{\vvv}{\boldsymbol{v}}
\newcommand{\eee}{\boldsymbol{e}}
\newcommand{\fff}{\boldsymbol{f}}
\newcommand{\www}{\boldsymbol{w}}
\newcommand{\0}{\boldsymbol{0}}
\newcommand{\lon}{\longrightarrow}
\newcommand{\ga}{\gamma}
\newcommand{\pa}{\partial}
\newcommand{\QED}{\hfill $\Box$}
\newcommand{\id}{{\mbox {\rm id}}}
\newcommand{\Ker}{{\mbox {\rm Ker}}}
\newcommand{\Image}{{\mbox {\rm Image}}}
\newcommand{\grad}{{\mbox {\rm grad}}}
\newcommand{\ind}{{\mbox {\rm ind}}}
\newcommand{\rot}{{\mbox {\rm rot}}}
\newcommand{\diver}{{\mbox {\rm div}}}
\newcommand{\Gr}{{\mbox {\rm Gr}}}
\newcommand{\GL}{{\mbox {\rm GL}}}
\newcommand{\LG}{{\mbox {\rm LG}}}
\newcommand{\Diff}{{\mbox {\rm Diff}}}
\newcommand{\Symp}{{\mbox {\rm Symp}}}
\newcommand{\Ct}{{\mbox {\rm Ct}}}
\newcommand{\Uns}{{\mbox {\rm Uns}}}
\newcommand{\rank}{{\mbox {\rm rank}}}
\newcommand{\sign}{{\mbox {\rm sign}}}
\newcommand{\Spin}{{\mbox {\rm Spin}}}
\newcommand{\Sp}{{\mbox {\rm Sp}}}
\newcommand{\Int}{{\mbox {\rm Int}}}
\newcommand{\Hom}{{\mbox {\rm Hom}}}
\newcommand{\Tan}{{\mbox {\rm Tan}}}
\newcommand{\codim}{{\mbox {\rm codim}}}
\newcommand{\ord}{{\mbox {\rm ord}}}
\newcommand{\Iso}{{\mbox {\rm Iso}}}
\newcommand{\corank}{{\mbox {\rm corank}}}
\def\mod{{\mbox {\rm mod}}}
\newcommand{\pt}{{\mbox {\rm pt}}}
\newcommand{\qed}{\hfill $\Box$ \par}
\newcommand{\spe}{\vspace{0.4truecm}}
\renewcommand{\0}{\mathbf 0}
\newcommand{\ad}{{\mbox{\rm ad}}}
\newcommand{\xdownarrow}[1]{%
  {\left\downarrow\vbox to #1{}\right.\kern-\nulldelimiterspace}
}

\newcommand{\dint}[2]{{\displaystyle\int}_{{\hspace{-1.9truemm}}{#1}}^{#2}}

\title{Topology of complements to real affine space line arrangements 
}

\author{Goo \textsc{Ishikawa}\thanks{Department of Mathematics, Faculty of Sciences, Hokkaido University, Sapporo 060-0810, Japan. 
e-mail: 
ishikawa@math.sci.hokudai.ac.jp}, 
\, 
Motoki \textsc{Oyama}\thanks{Department of Mathematics, Faculty of Sciences, Hokkaido University, Sapporo 060-0810, Japan. 
e-mail: oyama.motoki@gmail.com
}}


%

\renewcommand{\thefootnote}{\fnsymbol{footnote}}
\footnotetext{
\\
\qquad
The first author was supported by JSPS KAKENHI No.15H03615.}

\date{}

\maketitle

\begin{abstract} 
It is shown that the diffeomorphism type of the complement to a real space line arrangement 
in any dimensional affine ambient space is determined only by the number of lines 
and the data on multiple points.  
\end{abstract} 

\section{Introduction}

Let ${\mathcal A} = \{ \ell_1, \ell_2, \dots, \ell_d\}$ be a real space 
line arrangement, or a configuration, 
consisting of affine $d$-lines in $\R^3$. 
The different lines $\ell_i, \ell_j (i \not= j)$ may intersect, so that the 
union 
$\cup_{i=1}^d \ell_i$ 
is an affine real algebraic curve of degree $d$ in $\R^3$ possibly with multiple points. 
In this paper we determine the topological type of the complement 
$M({\mathcal A}) := \R^3 \setminus 
(\cup_{i=1}^d \ell_i)$ 
of ${\mathcal A}$, which is an open 
$3$-manifold. We observe that the topological type $M({\mathcal A})$ is determined only 
by the number of lines 
and the data on multiple points of ${\mathcal A}$. 
Moreover we determine the diffeomorphism type of $M({\mathcal A})$. 

Set $D^n := \{ x \in \R^n \mid \Vert x\Vert \leq 1\}$, the $n$-dimensional closed disk. 
The pair $(D^i \times D^j, D^i\times \pa(D^j))$
with $i + j = n, 0 \leq i, 0 \leq j$, 
is called an $n$-dimensional {\it handle of index} $j$ (see \cite{Wall}\cite{GS} for instance).  

\begin{figure}[ht]
\begin{center}
\includegraphics[width=5truecm, height=2.5truecm, clip, 
]{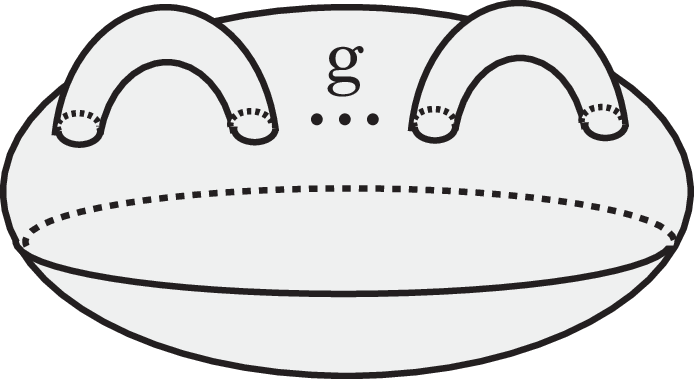}
\qquad 
\includegraphics[width=6truecm, height=2.5truecm, clip, 
]{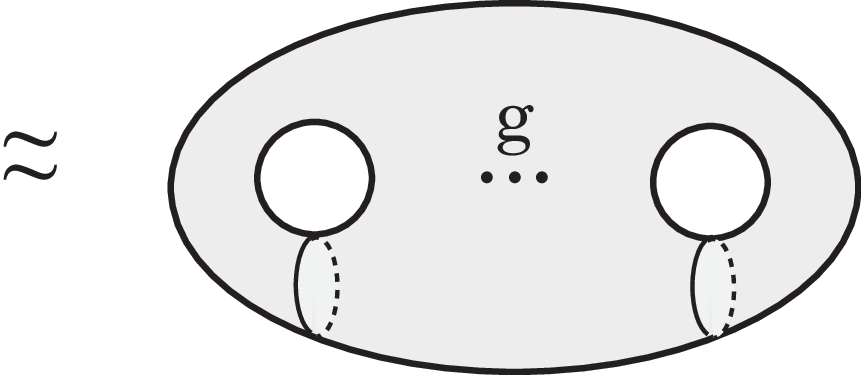}
\caption{$3$-ball with trivial $g$-handles of index $1$.}
\label{ball-with-handles}
\end{center}%
\end{figure}

Now take one $D^3$ 
and, for any non-negative integer $g$, 
attach to it $g$-number of $3$-dimensional handles 
$(D_k^2 \times D_k^1, D_k^2 \times \pa(D_k^1))$ of index $1$ ($1 \leq k \leq g$),  
by an attaching embedding $\varphi : \bigsqcup_{k=1}^g (D_k^2 \times \pa(D_k^1)) \to \pa(D^3) = S^2$ such that the obtained $3$-manifold 
$$
{\textstyle 
B_g := D^3 \bigcup_\varphi (\bigsqcup_{k=1}^g (D_k^2 \times D_k^1))}
$$ 
is orientable. 
We call $B_g$ 
the {\it $3$-ball with trivial $g$-handles of index $1$} (Figure \ref{ball-with-handles}.)

Note that the topological type of $B_g$ does not depend on the attaching map $\varphi$ and  
is uniquely determined only by the number $g$. The boundary of $B_g$ 
is the orientable closed surface $\Sigma_g$ of genus $g$. 

Let ${\mathcal A}$ be any $d$-line arrangement in $\R^3$. 
Let $t_i = t_i({\mathcal A})$ denote the number of multiple points with multiplicity $i$, 
$i = 2, \dots, d$. The vector $(t_d, t_{d-1}, \dots, t_2)$ provides a degree of 
degeneration of the line arrangement ${\mathcal A}$. 
Set 
$
g := d + \sum_{i=2}^d (i - 1)t_i. 
$
In this paper we show the following result:  

\bet
\label{three-dimensional-affine-topology}
The complement 
$M({\mathcal A})$ is homeomorphic to the interior of $3$-ball with trivial $g$-handles of index $1$. 
\ent

\bec
\label{three-dimensional-affine-homotopy}
$M({\mathcal A})$ is homotopy equivalent to the bouquet $\bigvee_{k=1}^g S^1$. 
\enc

The above results are naturally generalised to any line arrangements in $\R^n (n \geq 3)$. 

Let ${\mathcal A} = \{ \ell_1, \ell_2, \dots, \ell_d\}$ be a line arrangement in $\R^n$ and 
set $M({\mathcal A}) := \R^n \setminus (\cup_{i=1}^d \ell_i)$. 
Again let $t_i$ denote the number of multiple points of ${\mathcal A}$ of multiplicity $i$, $i = 2, \dots, d$. 
Set $g := d + \sum_{i=2}^d (i - 1)t_i$. 
Then we have 

\bet
\label{general-dimensional-affine-topology}
$M({\mathcal A})$ is homeomorphic to the interior of $n$-ball $B_g$ 
with trivially attached $g$-handles of index $n-2$. 
\ent

Thus we see that the topology of complements of real space line arrangements is 
completely determined by the combinational data, the {\it intersection poset} in particular. 
Recall that the intersection poset $P = P({\mathcal A})$ is the 
partially ordered set which consists of all multiple points, 
the lines themselves $\ell_1, \ell_2, \dots, \ell_d$ and $T = \R^n$ as elements, endowed with the inclusion order. 
Then the number $t_i$ is recovered as the number of minimal points $x$ such that $\#\{ y \in P \mid x < y, y \not= T\} = i$ and $d$ as the number of maximal points of $P \setminus \{ T\}$.

\bec
\label{general-dimensional-affine-homotopy}
$M({\mathcal A})$ is homotopy equivalent to the bouquet $\bigvee_{k=1}^g S^{n-2}$. 
\enc

In particular $M({\mathcal A})$ is a minimal space, i.e. 
it is homotopy equivalent to a $CW$ complex such that the number of $i$-cells is 
equal to its $i$-th Betti number for all $i \geq 0$. 

Even for semi-algebraic open subsets in $\R^n$, homotopical equivalence does not 
imply topological equivalence in general. 
However we see this is the case for complements of real affine line 
arrangements, as a result of Theorem 
\ref{general-dimensional-affine-topology} 
and Corollary \ref{general-dimensional-affine-homotopy}.

By the uniqueness of smoothing of corners,  
and by careful arguments at all steps of the proof of 
Theorem \ref{general-dimensional-affine-topology}, 
we see that 
Theorem \ref{general-dimensional-affine-topology} can be proved in differentiable category. 

\bet
\label{general-dimensional-affine-diffeomorphism}
$M({\mathcal A})$ is diffeomorphic to the interior of $n$-ball $B_g$ 
with trivially attached $g$-handles of index $n-2$. 
\ent

Note that the relative classification problem of line arrangements 
$(\R^n, \cup_{i=1}^d \ell_i)$ is classical but far from being solved (\cite{Hirschfeld} for instance). 
Moreover it has much difference in differentiable category and topological category. 
In fact even the local classification near multiple points of high multiplicity $i$, $i \geq n+2$
has moduli in differentiable category while it has no moduli in topological category. 
The classification of complements turns to be easier and simpler as we 
observe in this paper.

The real line arrangements on the plane $\R^2$ is one of classical and interesting 
subjects to study. 
It is known or easy to show that the number of connected components of the complement to 
a real plane line arrangement is given exactly by $1 + g$ using the number 
$g = d + \sum_{i=2}^d (i - 1)t_i$. This can be derived from 
Corollary \ref{general-dimensional-affine-homotopy} by just setting $n = 2$. 
For example, it can be shown from known combinatorial results for line arrangements on 
projective plane (see \cite{Grunbaum} for instance). In fact 
we prove it using our method in the process of the proof of 
Theorem \ref{general-dimensional-affine-topology}. 
Therefore Theorem \ref{general-dimensional-affine-topology} and 
Corollary \ref{general-dimensional-affine-homotopy} are regarded as 
a natural generalisation of the classical fact.

Though our object in this paper is the class of real affine line arrangements, 
it is natural to consider also real projective line arrangements 
consisting of projective lines in the projective space $\R P^n$, 
or corresponding real linear plane arrangements 
consisting of $2$-dimensional linear subspaces 
in $\R^{n+1}$. However the topology of complements 
in both cases are not determined, in general, by 
the intersection posets, which are defined similarly to the affine case. 
In fact it 
is known an example of pairwise transversal linear plane arrangements 
${\mathcal B}$ and ${\mathcal B'}$ in $\R^4$ with $d = 4$ 
such that the complements 
$M({\mathcal B})$ and $M({\mathcal B'})$ have non-isomorphic cohomology algebras 
and therefore they are not homotopy equivalent, so, not homeomorphic to each other 
(\cite{Ziegler}, Theorem 2.1). 

A linear plane arrangement in $\R^4$ is pairwise transverse if 
and only if the corresponding projective line arrangement is 
non-singular (without multiple points) in $\R P^3$. 
Non-singular line arrangements in $\R P^3$, which are called skew line configurations, 
are studied in details (see \cite{Hirschfeld, Viro85, VD, VV} for instance). 
Moreover, the topology of non-singular real algebraic curves in $\R P^3$ is studied, 
related to Hilbert's 16th problem, by many authors (see \cite{MO} for instance). 
Also refer to the surveys on the study of real algebraic varieties (\cite{Gudkov, Viro86}). 

It is natural to consider also complex line arrangements in $\C^n = \R^{2n}$. 
The topology of complex subspace arrangements in $\C^n$, 
in particular, homotopy types of them is studied in 
detail (see \cite{OT, Ziegler} for instance). 
Then it is known that the intersection poset turns to have much information 
in complex cases than in real cases. 
Refer to \cite{Vasiliev, ZZ}, for instance, on the theory on the homotopy types of complements 
for general subspace arrangements.

In \S \ref{Trivial handle attachments},
we define the notion of trivial handle attachments clearly. 
In \S \ref{Affine line arrangements}, we show Theorem \ref{general-dimensional-affine-topology} 
and Theorem \ref{general-dimensional-affine-diffeomorphism} in parallel, 
using an idea of stratified Morse theory (\cite{GM}) in a simple situation. 
We then realize a deference of topological features between the 
complements to line arrangements and to knots, links, tangles or general spacial graphs 
(Remark \ref{knot-complement}). 
In the last section, related to our results, 
we discuss briefly the topology of real projective line arrangements and 
real linear plane arrangements.

The authors thank Masahiko Yoshinaga for his valuable suggestion to 
turn authors' attention to real space line arrangements. They thank also an 
anonymous referee for his/her valuable comments.

\section{Trivial handle attachments}
\label{Trivial handle attachments}

First we introduce the local model of trivial handle attachments. 

Let $j < n$. Let $S^j \subset \R^n$ be the sphere defined by $x_1^2 + \cdots + x_j^2 + x_n^2 = 1, 
x_{j+1} = 0, \dots, x_{n-1} = 0$, and $\pa(D^j) = S^{j-1} = S^j \cap \{ x_n = 0\}$. 
Let $e_\ell \in \R^n$ be the vector defined by $(e_\ell)_i = \delta_{\ell i}$. 
Then define an embedding $\widetilde{\Phi} : 
D^{n-j}\times S^j \to \R^n$ by 
$$
\widetilde{\Phi}(t_1, \dots, t_{n-j-1}, t_{n-j}, x) 
:= x + t_1e_{n-1} + \cdots + t_{n-j-1}e_{j+1} + t_{n-j}x, 
$$
which gives a tubular neighbourhood of $S^j$ in $\R^n$. Set 
$$
\varphi_{\mbox{\rm{st}}} := \widetilde{\Phi}\vert_{D^{n-j}\times\pa(D^j)} : D^{n-j}\times S^{j-1} 
\to \R^{n-1} \subset \R^n, 
$$
which gives a tubular neighbourhood of $S^{j-1}$ in $\R^{n-1} = \{ x_n = 0\}$.
We call $\varphi_{\mbox{\rm{st}}}$ the {\it standard} attaching map of the handle of index $j$. 
Note that the embedding $\varphi_{\mbox{\rm{st}}}$ extends to the {\it standard} handle 
$\Phi : D^{n-j}\times D^j \to \R^n$, which is defined by 
$$
\Phi(t_1, \dots, t_{n-j-1}, t_{n-j}, x_1, \dots, x_j) := 
\widetilde{\Phi}\left(t_1, \dots, t_{n-j-1}, t_{n-j}, x_1, \dots, x_j, 0, \dots, 0, 
\textstyle{\sqrt{1 - \sum_{i=1}^j x_i^2}}\right), 
$$
attached to $\{ x_n \leq 0 \}$ along $\varphi_{\mbox{\rm{st}}}$. 

Let $M$ be a topological (resp. differentiable) 
$n$-manifold with a {\it connected} boundary $\pa M$.  

Let $p \in \pa M$. 
A coordinate neighbourhood $(U, \psi)$, $\psi : 
U \to \psi(U) \subset \R^{n-1}\times\R$ around $p$ in $M$ 
is called {\it adapted} if $\psi : U \to \R^n$ is a homeomorphism of $U$ 
and $\psi(U) \cap \{ x_n \leq 0\}$ which maps $U \cap \pa M$ to 
$\R^{n-1} = \{ x_n = 0\}$. 

\

\begin{figure}[ht]
\begin{center}
\includegraphics[width=12truecm, height=2.5truecm, clip, 
]{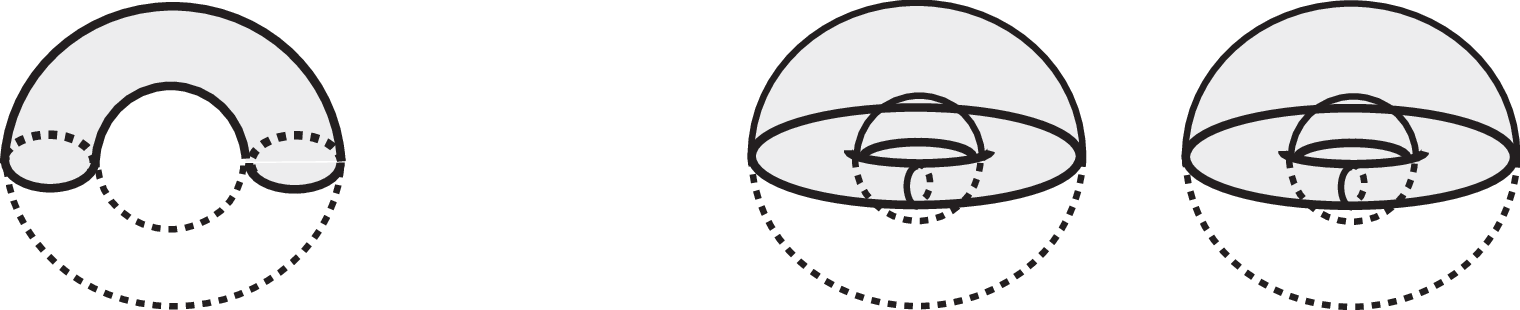} 
\caption{Trivial handle attachments: the cases $n = 3, j = 1, \ell = 1$ and $n = 4, j = 2, \ell = 2$.}
\label{Unknotted}
\end{center}%
\end{figure}

Now we consider an attaching of several number of handles of index $j$ to $M$ along $\pa M$. 
We call a handle attaching map $\varphi : \bigsqcup_{k=1}^\ell (D^{n-j}_k\times \pa(D^j_k))
\to \pa M$ {\it trivial} if there exist disjoint adapted coordinate neighbourhoods 
$(U_1, \psi_1), \dots, (U_\ell, \psi_\ell)$ on $M$ such that 
$\varphi(D^{n-j}_k\times \pa(D^j_k)) \subset U_k$ and 
$\psi_k\circ\varphi : D^{n-j}_k\times \pa(D^j_k) \to \R^{n-1}\times \R$ is the standard 
attachment for $k = 1, \dots, \ell$. 
 (Figure \ref{Unknotted})

Then $M \cup_{\varphi} \left(\bigsqcup_{k=1}^\ell (D^{n-j}_k\times D^j_k) \right)$ 
is called the manifold obtained from $M$ by attaching standard handles 
and the topological type of $M$ does not depends on the attaching map $\varphi$ but 
depends only on $j$ and $\ell$. 
Moreover if $M$ is a differentiable manifold, the diffeomorphism type 
of the attached manifold is uniquely determined by the smoothing or straightening of corners
(see Proposition 2.6.2 of \cite{Wall} for instance). 
Note that the diffeomorphism type of the interior does not change by the smoothing. 

\begin{figure}[ht]
\begin{center}
\includegraphics[width=15truecm, height=1.8truecm, clip, 
]{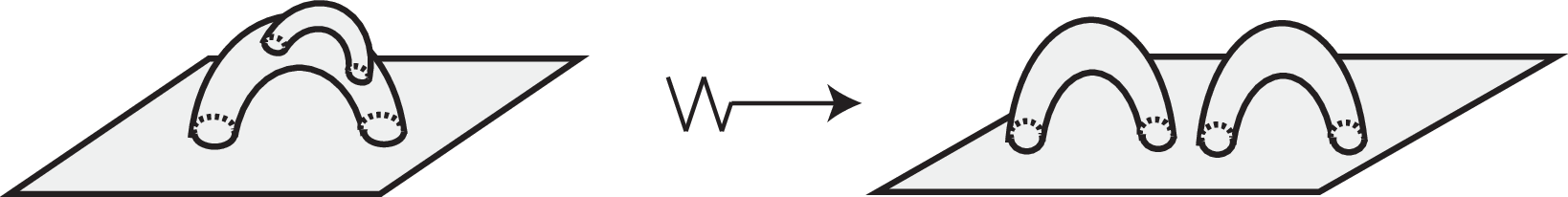} 
\caption{Sliding of trivial handle attachments.}
\label{Handle-attach}
\end{center}%
\end{figure}

Note that, if $\varphi$ is a trivial handle attaching map, then $\varphi\vert_{0\times\pa(D^j_k)} 
: 0\times\pa(D^j_k) \to \pa M$ is unknotted and 
$\varphi\vert_{\bigsqcup_{k=1}^\ell (0\times\pa(D^j_k))} 
: \bigsqcup_{k=1}^\ell (0\times\pa(D^j_k)) \to \pa M$ is unlinked (see Figure \ref{Handle-attach}). 
Therefore we can slide the trivial attachment mapping $
\bigsqcup_{k=1}^\ell (D^{n-j}_k\times \pa(D^j_k))$ to an 
embedding into a disjoint union to an arbitrarily small neighbourhoods of any disjoint $\ell$ number points on $\pa M$ up to isotopy (cf. Homogeneity Lemma \cite{Milnor}). 

\ber
{\rm 
The assumption that $\pa M$ is connected is essential. For example, let $M = \{ x \in \R^n \mid 
-1 \leq x_n \leq 1\}$. 
Then we have at least two non-homeomorphic spaces by different 
attachments of two trivial handles of index $1$ (Figure \ref{non-connected-boundary}). 

\begin{figure}[ht]
\begin{center}
\includegraphics[width=12truecm, height=2.5truecm, clip, 
]{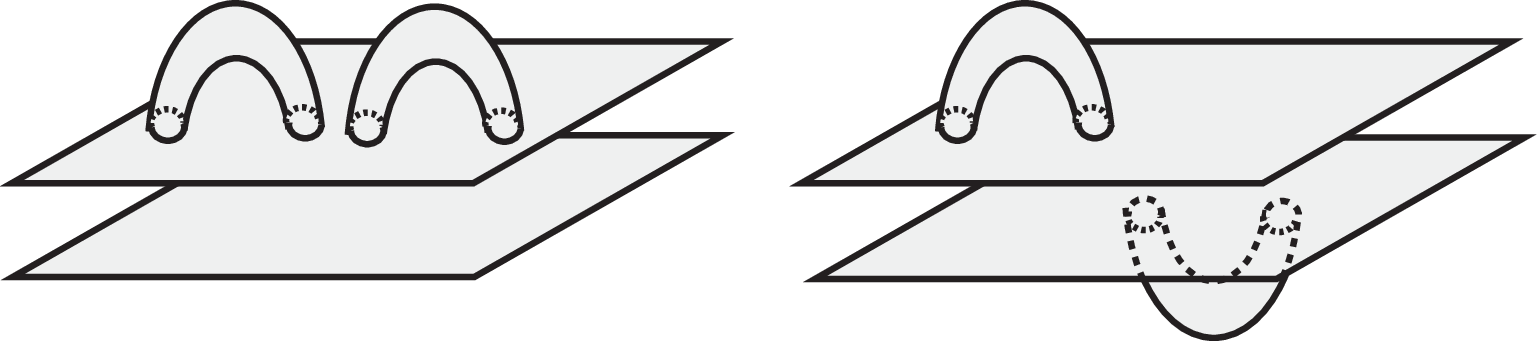} 
\caption{Non-homeomorphic attachments of trivial handles $n = 3, j = 1, \ell = 2$.}
\label{non-connected-boundary}
\end{center}%
\end{figure}%

}
\enr

We see that iterative trivial attachments result a homeomorphic (resp. differentiable) 
manifold to a simultaneous trivial attachments. 

\bel
\label{iterative-attachments}
Let $M'$ be a topological (resp. differentiable) 
$n$-manifold with connected boundary $\pa M'$. 
Suppose $M'$ is homeomorphic (diffeomorphic) to a space 
$M_1 := M \cup_{\varphi} \left(\bigsqcup_{k=1}^\ell (D^{n-j}_k\times D^j_k) \right)$
obtained, from a topological (differentiable) manifold $M$ with connected boundary, 
by attaching $k$ number of trivial handles of index $j$. 
Then the space $M_2 := M' \cup_{\varphi'} \left(\bigsqcup_{k=\ell+1}^{\ell+m} 
(D^{n-j}_k\times D^j_k) \right)$ obtained from $M'$ by attaching 
$m$ number of trivial handles of index $j$ is homeomorphic (diffeomorphic) to the space 
$M_3 := M \cup_{\varphi''} \left(\bigsqcup_{k=1}^{\ell+m} (D^{n-j}_k\times D^j_k) \right)$ 
obtained from $M$ by attaching $\ell + m$ number of 
trivial handles of index $j$. 
\enl

See Figure \ref{Handle-attach} for the case $j = 1$. 

\

\noindent
{\it Proof of Lemma \ref{iterative-attachments}.}
Let $f : M_1 \to M'$ be a homeomorphism (resp. a diffeomorphism). 
Then $f(\bigsqcup_{k=1}^\ell (D^{n-j}_k\times D^j_k)$ does not contained in $\pa M'$. 
Then we slide, up to isotopy, the attaching map $\varphi' : 
\bigsqcup_{k=\ell+1}^{\ell+m} 
(D^{n-j}_k\times \pa D^j_k) \to \pa M'$ to $\varphi''' : 
\bigsqcup_{k=\ell+1}^{\ell+m} 
(D^{n-j}_k\times \pa D^j_k) \to \pa M'$ such that 
$$
f\left( \varphi\left( \bigsqcup_{k=1}^\ell \left( D^{n-j}_k\times \pa D^j_k\right)\right)\right) \cap \varphi'''\left(\bigsqcup_{k=\ell+1}^{\ell+m} 
\left( D^{n-j}_k\times \pa D^j_k\right)\right)
= \emptyset. 
$$
Consider $\varphi'' := \varphi \bigsqcup f^{-1}\circ \varphi''' 
: \bigsqcup_{k=1}^{\ell+m} (D^{n-j}_k\times \pa D^j_k) \to \pa M$. 
Then $M_2$ is homeomorphic (resp. diffeomorphic) to $M_3$. 
\QED

\section{Affine line arrangements}
\label{Affine line arrangements}

Let $n \geq 2$. 

\begin{figure}[ht]
\begin{center}
\includegraphics[width=5truecm, height=2.5truecm, clip, 
]{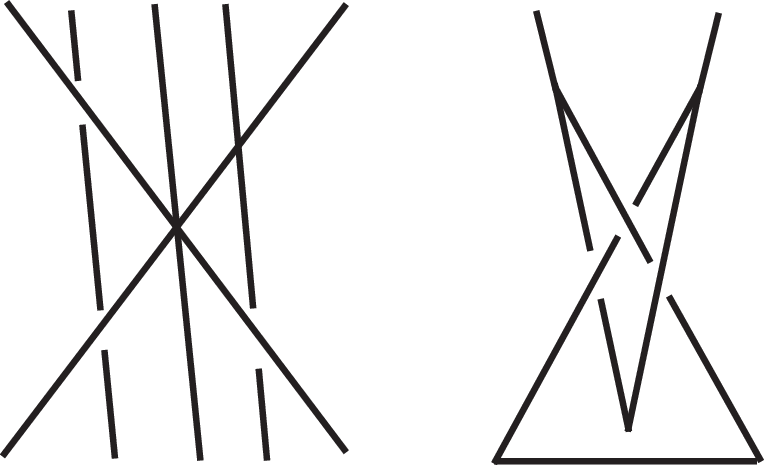} 
\caption{A line arrangement and a space graph}
\label{space-graph}
\end{center}
\end{figure}

We consider line arrangements in $\R^n$ or more generally consider 
a subset $X$ in $\R^n$ which is a union of finite number of closed 
line segments and half lines. 
Then $X$ may be regarded as a finite graph (with non-compact edges)  
embedded as a closed set in $\R^n$ (Figure \ref{space-graph}). Here we admit vertices of valency $1$. 


Take a unit vector $v \in S^{n-1} \subset \R^n$ and define the height function 
$h : \R^n \to \R$ by $h(x) := x\cdot v$ using Euclidean inner product. 
Choose $v$ so that 

(i) $v$ is neither perpendicular to any line segments nor half lines in $X$. 

(ii) For each $c$, the hyperplane $h(x) = c$ of level $c$ contains at most one vertex of $X$. 

Note that there exists a union $\Sigma$ of finite number of great hyperplanes such that 
any unit vector in $S^{n-1} \setminus \Sigma$ satisfies the conditions (i) and (ii). 

After a rotation of $\R^n$, we may suppose $h(x) = x_n$. We write $x = (x', x_n)$, 
where $x' = (x_1, \dots, x_{n-1})$. 
Set $M = \R^n \setminus X$ and, for any $c \in \R$, 
$$
M_{\leq c} := \{ x \in M \mid x_n \leq c\}, \quad 
M_{< c} := \{ x \in M \mid x_n < c\}.
$$
Let $V \subset X$ be the set of vertices of $X$. Set 
$V = \{ u_1, u_2, \dots, u_r\}, c_i = h(u_i)$ and $C = h(V) = \{ c_1, c_2, \dots, c_r\}$ 
with $c_1 < c_2 < \cdots < c_r$. 

\

Though the following lemma 
is clear intuitively, we give a proof to make sure. 

\bel
\label{Topological-triviality}
The topological (resp. diffeomorphism) type of $M_{\leq c}$ is constant on $c_i < c < c_{i+1}$  
and 
the topological (diffeomorphism) 
type of $M_{< c}$ is constant on $c_i < c \leq c_{i+1}$, $i = 0, 1, \dots, r$, 
with $c_0 = -\infty, c_{r+1} = \infty$. Here $M_{< \infty}$ means $M$ itself. 
\enl

\Proof
First we treat the case $i < r$. 
Take a sufficiently large $R > 0$ such that 
$\{ x \in X \mid c_i < x_n < c_{i+1}, \Vert x'\Vert > R/2 \} = \emptyset$. 
Consider the cylinder 
$$
C := \{ x \in \R^n \mid c_i < x_n < c_{i+1}, \Vert x'\Vert \leq R\}. 
$$ 
Then ${\mathcal C} := \{ {\mbox{\rm Int}}\, C \setminus X, X\cap C, \pa C \}$ 
is a Whitney stratification of $C$. 
The function $h : C \to (c_i, c_{i+1})$ is proper and the restriction of $h$ to each stratum is 
a submersion. Now we follow the standard method (the proof of Thom's first isotopy lemma 
\cite{Thom, Mather}) to show differentiable triviality of mappings. 
Note that the flow used in the proof of isotopy lemma is differentiable in each stratum. 
For any $\varepsilon > 0$, 
take a vector field $\eta$ over $(c_i, c_{i+1})$ such that $\eta = 0$ on 
$(c_i, c_i + \varepsilon/2)$ and $\eta = \pa/\pa y$ on $(c_i + \varepsilon, c_{i+1})$, 
where $y$ is the coordinate on $\R$. 
Then $\eta$ lifts to a controlled vector field $\xi$ over $C$ such that $\xi$ tangents to each stratum. We extend $\xi\vert_{\pa c}$ to $\{ x \in \R^n \mid c_i < x_n < c_{i+1}, \Vert x'\Vert \geq R\}$ 
via the retraction 
$x = (x', x_n) \mapsto (\frac{1}{\Vert x'\Vert}Rx', x_n)$ and 
to $\{ x \in \R^n \mid x_n < c_i + \varepsilon/2\}$ by letting it $0$, and 
we have an integrable vector field $\xi$ on 
$\{ x \in \R^n \mid x_n < c_{i+1} \}$. By integrating $\xi$, we have a homeomorphism 
of $M_{\leq c}$ and $M_{\leq c'}$ for any $c, c' \in (c_i, c_{i+1})$ and 
a diffeomorphism 
of $M_{< c}$ and $M_{< c'}$ for any $c, c' \in (c_i, c_{i+1}]$. Note that the differentiable flow 
of the vector field may not be defined through $x_n = c_{i+1}$ but it gives 
a diffeomorphism of  $M_{< c}$ and $M_{< c_{i+1}}$. 

Second we treat the case $i = r$. 
Consider the quadratic cone $\Vert x'\Vert^2 - Rx_n^2 = 0$ in $\R^n$. 
Supposing $c_{r+1} > 0$ after a translation along $x_n$-axis in necessary, and taking $R$ sufficiently large, we have $X \cap \{ x \in \R^n \mid c_{r+1} < x_n\}$ lies inside of 
the cone $\Vert x'\Vert^2 - Rx_n^2 < 0$. 
Now set 
$$
D := \{ x \in \R^n \mid c_{r+1} < x_n, \Vert x'\Vert^2 - Rx_n^2 \leq 0 \}, 
$$
and consider the proper map $h : D \to (c_{r+1}, \infty)$ with the 
Whitney stratification 
${\mathcal D} := \{ {\mbox{\rm Int}}\, D \setminus X, X \cap D, \pa D \}$. 
For any $\varepsilon > 0$, 
take a (non-complete) vector field $\eta$ over $(c_{r+1}, \infty)$ such that $\eta = 0$ on 
$(c_{r+1}, c_{r+1} + \varepsilon/2)$ and 
$\eta = (1 + y^2)\pa/\pa y$ on $(c_{r+1}, \infty)$. 
We lift $\eta$ to a controlled vector filed $\xi$ over $D$ and then over $\R^n$. 
Then, using the integration of $\xi$, we have a diffeomorphism of 
$M_{\leq c}$ and $M_{\leq c'}$ for any $c, c' \in (c_i, c_{i+1})$, 
and a diffeomorphism of 
$M_{< c}$ and $M_{< c'}$ for any $c, c' \in (c_i, c_{i+1}]$. 
In particular 
we have that $M_{< c}$ for $c_{r+1} < c$ is diffeomorphic to $M$ itself. 
\QED

\

\ber
{\rm 
The topological (resp. diffeomorphism) 
type of $M_{\leq c}$ (resp. $h^{-1}(c) \setminus X$) is not necessarily constant at $c = c_{i+1}$. 
}
\enr

We observe the topological change of $M_{< c}$ when $c$ moves across a critical value $c_i$ 
as follows:

\bel
\label{topological-bifurcation}
Let $u$ be a vertex of $X$ and let $c = h(u)$. 
Let $s = s(u)$ denote 
the number of edges of $X$ which are adjacent to $u$ from above 
with respect to $h$. 

Then, 
for a sufficiently small $\varepsilon > 0$, 
the open set $M_{< c+\varepsilon}$ is diffeomorphic to the interior of 
$M_{\leq c - \varepsilon} \bigcup_\varphi (\bigsqcup_{i=1}^{s-1} (D_i^2 \times D_i^{n-2}))$, obtained 
by an attaching map 
$$
\varphi : \bigsqcup_{i=1}^{s-1} D^2 \times \pa(D^{n-2}) \longrightarrow 
h^{-1}(c - \varepsilon) \setminus X\ = \pa(M_{\leq c - \varepsilon}) \subset M_{\leq c - \varepsilon}, 
$$ 
of $(s-1)$ number of trivial handles of index $n-2$, provided $s \geq 1$. 

In particular 
$M_{< c+\varepsilon}$ is diffeomorphic to $M_{< c - \varepsilon}$ if $s = 1$. 

\begin{figure}[ht]
\begin{center}
\includegraphics[width=12truecm, height=2.5truecm, clip, 
]{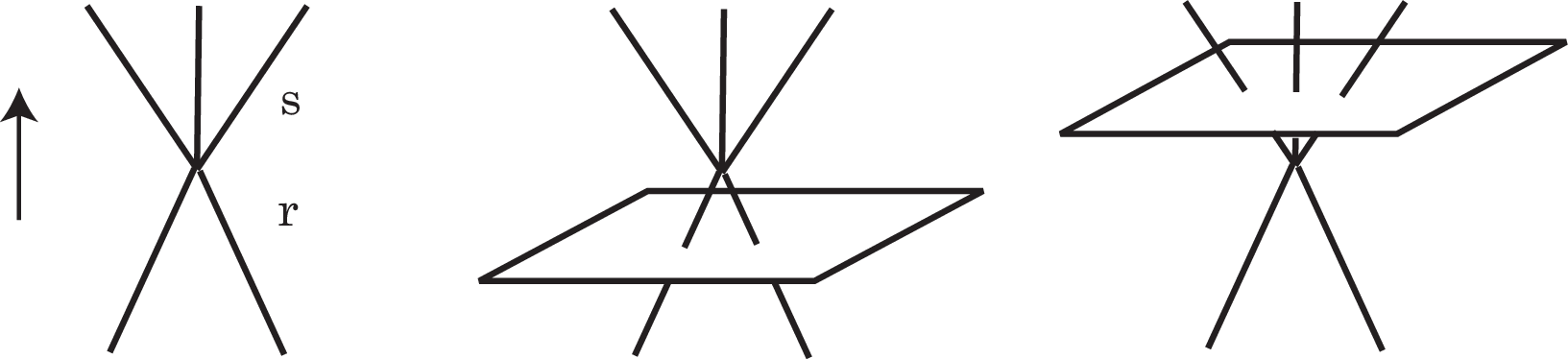} 
\caption{Topological bifurcations.}
\label{Topological-bifurcation}
\end{center}%
\end{figure}

If $s = 0$ then 
$M_{< c+\varepsilon}$ is diffeomorphic to the interior of 
$M_{\leq c - \varepsilon} \bigcup_\varphi (D^1 \times D^{n-1})$ obtained 
by an attaching map $\varphi : D^1 \times \pa(D^{n-1}) \to 
h^{-1}(c - \varepsilon) \setminus X$ of a (not necessarily trivial) handle of index $n-1$. 
{\rm (See Figure \ref{Topological-bifurcation}. )}
\enl

\begin{figure}[ht]
\begin{center}
\includegraphics[width=7truecm, height=2truecm, clip, 
]{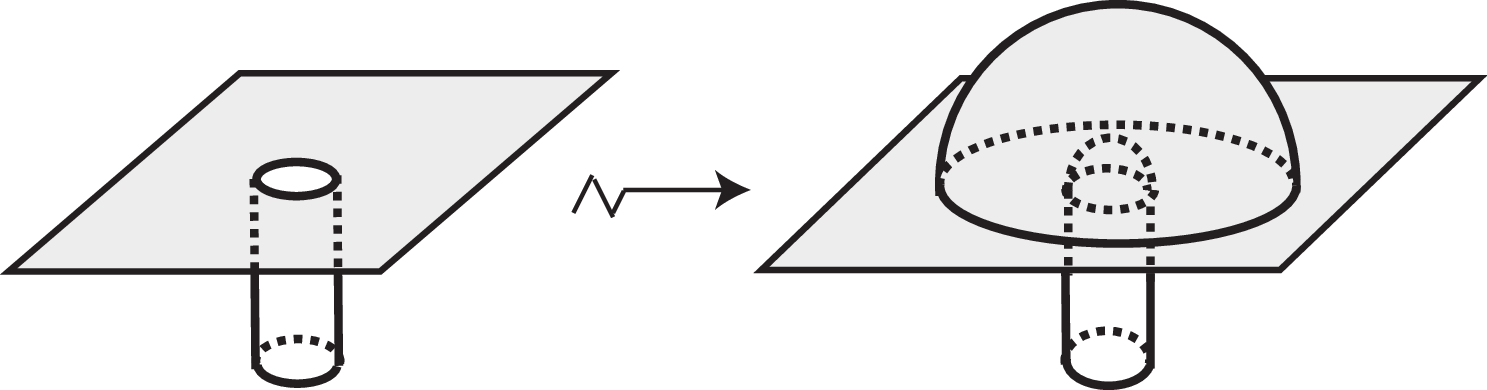} 
\caption{Topological change in the case $s = 0$.}
\label{change-5}
\end{center}
\end{figure}

\ber
{\rm 
In the case $s = 0$, the handle attachment is not necessarily trivial 
since the core of the attachment does not necessarily bounds a disk. (See Figure \ref{change-5}.) 
}
\enr

\ber
{\rm 
Note that if $r = r(u)$ denotes
the number of edges of $X$ which are adjacent to $p$ from {\it below} 
with respect to $h$, then 
the intersection $X \cap h^{-1}(c - \varepsilon)$ consists of $r$-points in 
the hyperplane $h^{-1}(c - \varepsilon)$ and thus 
$h^{-1}(c - \varepsilon) \setminus X$ is a punctured hyperplane by $r$-points. 
}
\enr

\ber
{\rm
Note that locally in a neighbourhood of each vertex $u$ of $X$, the topological equivalence class of 
the germ of a generic height function $h : (\R^n, X, u) \to (\R, c)$ 
is determined only by $s$ and $r$, the numbers of branches. This can be shown 
by using Thom's isotopy lemma (\cite{Mather}). 
}
\enr

\begin{figure}[ht]
\begin{center}
\includegraphics[width=13truecm, height=2.5truecm, clip, 
]{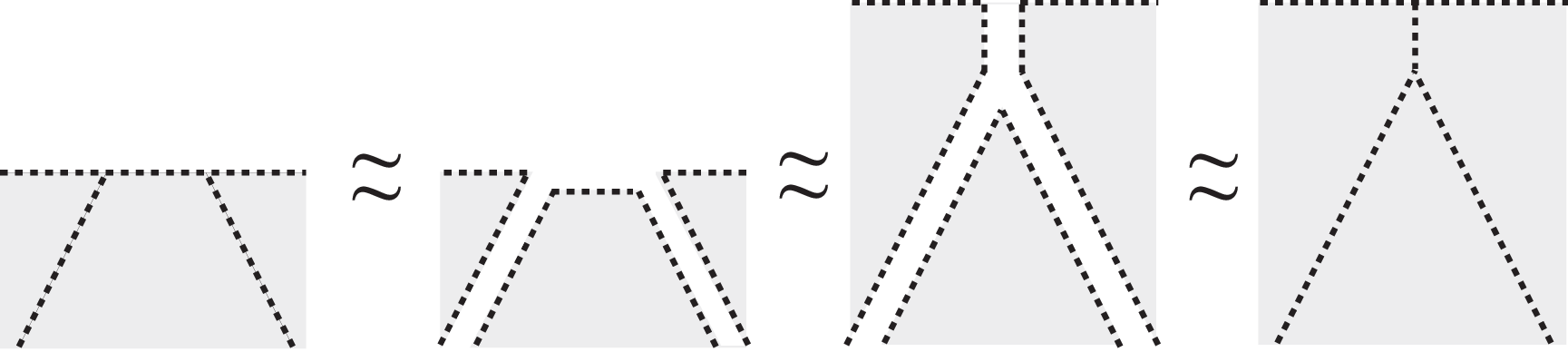} 
\caption{No topological changes of complements occur when $s = 1$.}
\label{no-change}
\end{center}
\end{figure}

\begin{figure}[ht]
\begin{center}
\includegraphics[width=11truecm, height=2.5truecm, clip, 
]{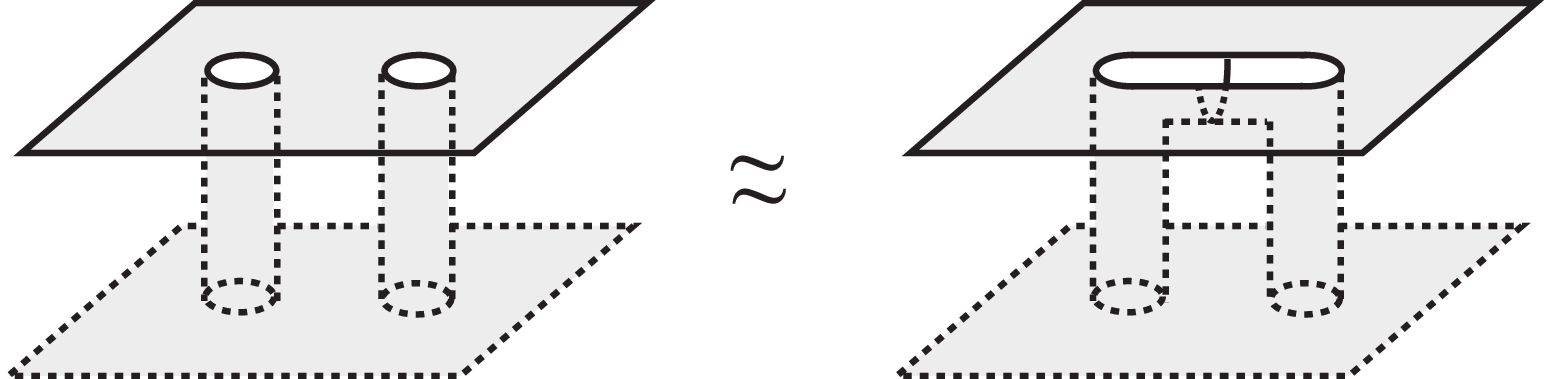} 
\caption{Boring a canal does not change the topology of ground.}
\label{no-change-2}
\end{center}
\end{figure}

\noindent
{\it 
Proof of Lemma \ref{topological-bifurcation}.}
For sufficiently small $0 < \varepsilon < \varepsilon'$, 
$M_{< c - \varepsilon} \setminus M_{\leq c - \varepsilon'}$ is a space 
$\{ x \in \R^n \mid c - \varepsilon' < h(x) < c - \varepsilon\}$ deleted $r$-half-lines. 
We may suppose the intersection $X \cap h^{-1}(c - \varepsilon)$ lies on a line, up to a diffeomorphism 
of $M_{\leq c - \varepsilon}$. 
We delete $r$-small tubular neighbourhoods of the half-lines from the half space, then still we have a diffeomorphic space to $M_{< c - \varepsilon} \setminus M_{\leq c - \varepsilon'}$. 
Then we connect the $r$-holes by boring a sequence of canals without changing the diffeomorphism type of complements. 
See Figures \ref{no-change} and \ref{no-change-2}. 
The boring a canal means, in general dimension, 
to delete $D^1 \times D^{n-1}$ along the line segment connecting the holes. 

First let $s = 1$. 
Then the resulting space 
is diffeomorphic to $M_{< c + \varepsilon} \setminus M_{\leq c - \varepsilon'}$. 
The diffeomorphism is taken to be the identity on  $M_{\leq c - \varepsilon'}$ and 
it extends to a diffeomorphism between $M_{< c - \varepsilon}$ and $M_{< c + \varepsilon}$. 
This shows Lemma \ref{topological-bifurcation} in the case $s = 1$. 

\begin{figure}[ht]
\begin{center}
\includegraphics[width=10truecm, height=2truecm, clip, 
]{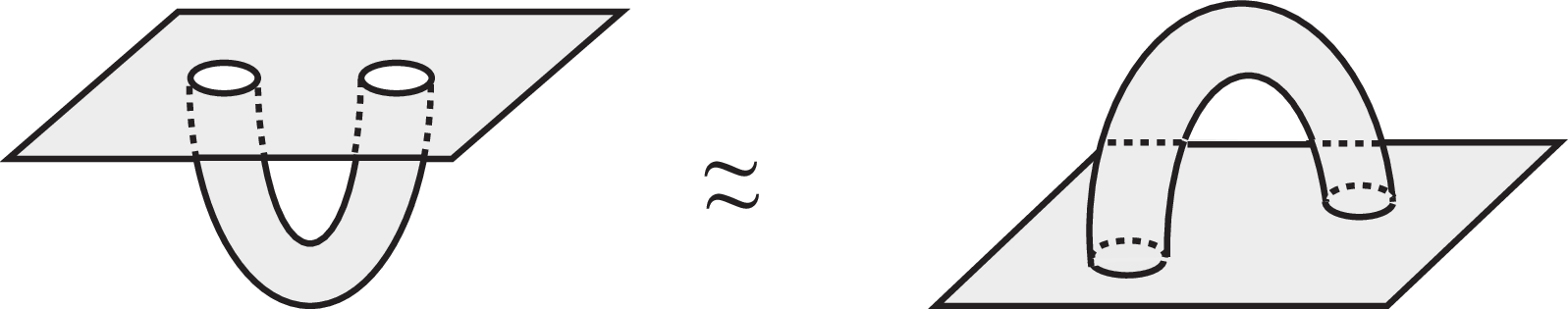} 
\caption{Digging a tunnel is same as bridging for the topology of ground.}
\label{change-1}
\end{center}
\end{figure}

\begin{figure}[ht]
\begin{center}
\includegraphics[width=10truecm, height=2.5truecm, clip, 
]{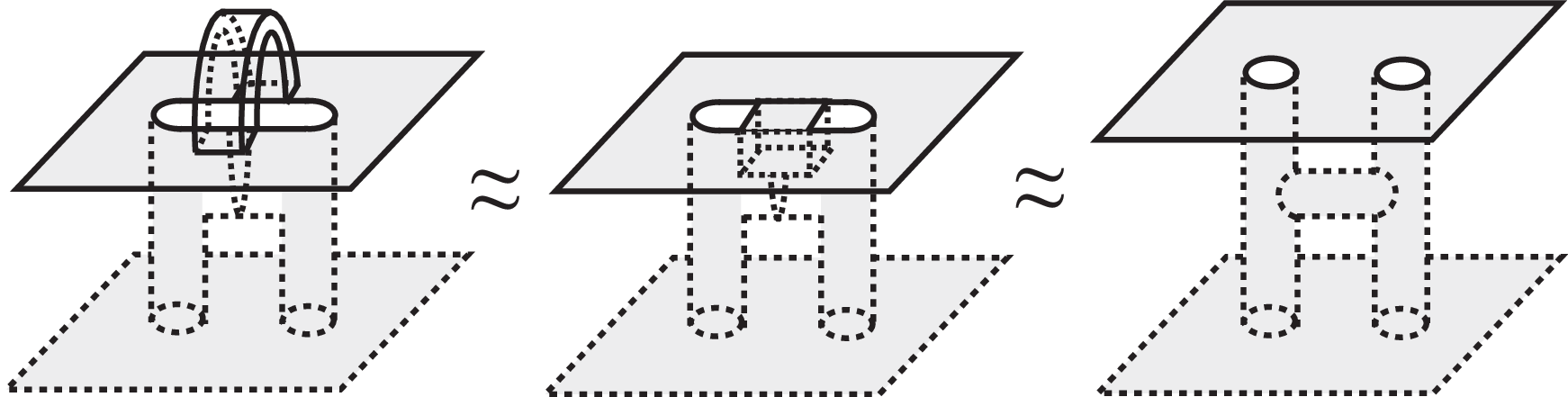} 
\caption{The case $s = 2, r = 2$.}
\label{change-2}
\end{center}
\end{figure}

Next we teat the case $s = 2, r = 0$. 
The topological change from $M_{c - \varepsilon}$ to $M_{c + \varepsilon}$ is give by 
digging a tunnel, which is, equivalently, given by a handle attaching of index $n-2$. 
In fact, we examine the topological change of the complement to 
$$
\sqcup = \{ (0, x_{n-1}, x_n) \in \R^n \mid (-2 \leq x_{n-1} \leq 2, x_n = 0) {\mbox{\rm \ or\ }} (x_{n-1} =-2, x_n \geq 0) {\mbox{\rm \ or\ }} (x_{n-1} = 2, x_n \geq 0)\}, 
$$
in $\R^n$ when across $x_n = c = 0$. 
Take the closed tube $T$ of radius $1$ of $\sqcup$. Then for the complement 
$M = \R^n \setminus T$,  
$M_{< \varepsilon}$ is diffeomorphic to the interior of the half space $\{ x_n \leq 0\}$ attached 
the handle
$$
H = \{ x \in \R^n \mid 
-1 \leq x_{n-1} \leq 1, \ \frac{1}{2} \leq x_1^2 + \cdots + x_{n-2}^2 + x_n^2 \leq 2, \ x_n \geq 0\}. 
$$
along 
$$
H \cap \{ x_n \leq 0\} = \{ x \in \R^n \mid 
-1 \leq x_{n-1} \leq 1, \ \frac{1}{2} \leq x_1^2 + \cdots + x_{n-2}^2 \leq 2\}. 
$$
The pair $(H, H \cap \{ x_n \leq 0\})$ 
is diffeomorphic to $(D^2\times D^{n-2}, D^2\times \pa D^{n-2})$ where 
the core $(0\times D^{n-2}, \pa D^{n-2})$ corresponds to 
$\{ x_1^2 + \cdots + x_{n-2}^2 + x_n^2 = 1, x_{n-1} = 0, x_n \geq 0\}$ and 
$\{ x_1^2 + \cdots + x_{n-2}^2 = 1, x_{n-1} = 0, x_n = 0\}$. 
Note that the latter bounds a $n-1$-dimensional disk 
$\{ x_1^2 + \cdots + x_{n-2}^2 \leq 1, x_{n-1} = 0, x_n = 0\}$, 
which does not touch the boundary $\pa M_{< \varepsilon}$. 
See Figures \ref{change-1} and \ref{change-2}. 

\begin{figure}[ht]
\begin{center}
\includegraphics[width=9truecm, height=2truecm, clip, 
]{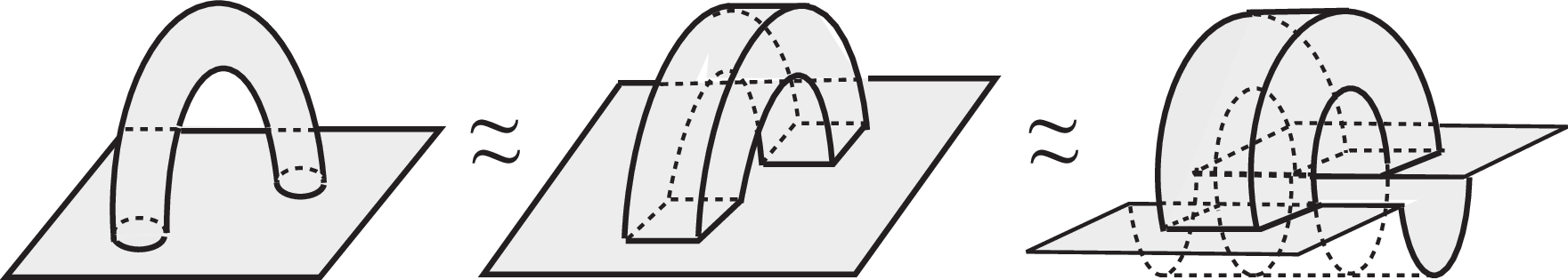} 
\caption{Trivial handle attachment and topological bifurcation.}
\label{change-3}
\end{center}
\end{figure}

\begin{figure}[ht]
\begin{center}
\includegraphics[width=10truecm, height=2.5truecm, clip, 
]{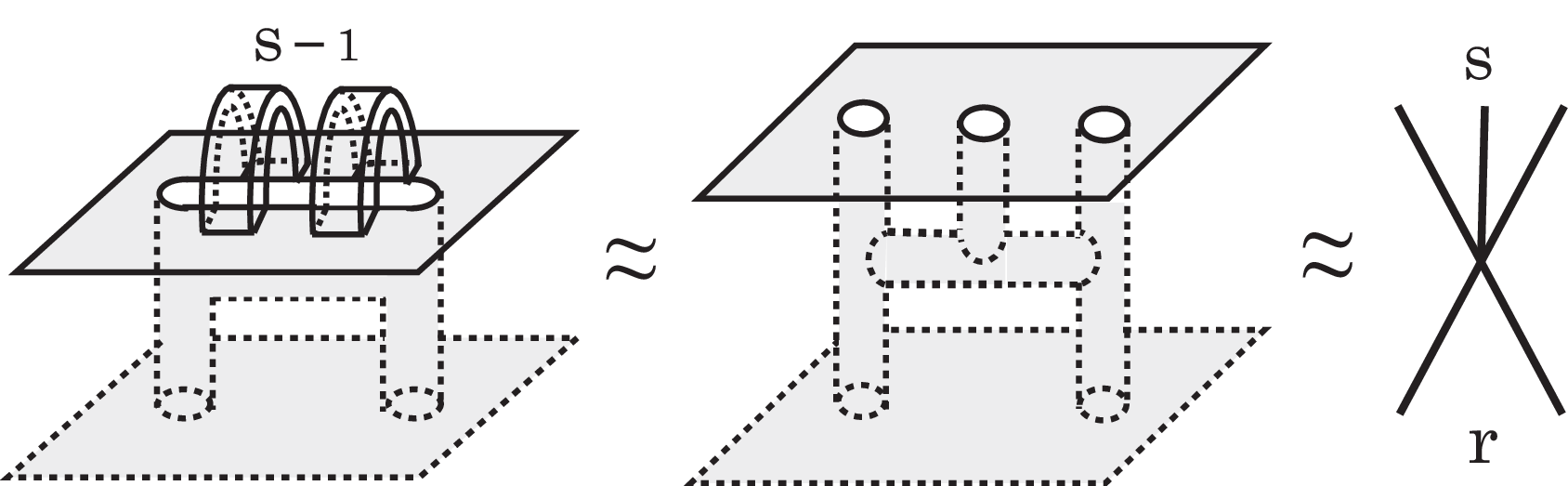} 
\caption{The case $s = 3, r = 2$.}
\label{change-4}
\end{center}
\end{figure}

The same argument works for any $r$. See Figure \ref{change-2} for the case $s = 2, r = 2$. 
Note that complements to \lq\lq X" and \lq\lq H" are diffeomorphic. 
See Figures \ref{change-2}, \ref{change-3} and \ref{change-4}. 


In general, for any $s \geq 2$, the topological change is obtained 
by attaching trivial $s-1$ handles of index $n - 2$. See Figure \ref{change-4}. 

In the case $s = 0$, contrarily to above, the change of diffeomorphism type is obtained by 
an attaching not necessarily trivial handle. See Figure \ref{change-5}. 

When $n = 2$, the topological bifurcation occurs just as putting 
$s-1$ number of disjoint open disks. 

Thus we have Lemma \ref{topological-bifurcation}. 
\QED

\

First let us apply Lemma \ref{Topological-triviality} and Lemma \ref{topological-bifurcation} 
to the case $n = 2$. 

For a $c \in \R$ of sufficiently large $\vert c\vert$, supposing a generic height function 
is given by $h = x_2$ as above. Then 
$M_{\leq c}$ (resp. $M_{< c}$) is diffeomorphic to the half plane $\{ x_2 \leq c \}$ 
(resp. $\{ x_n < c \}$ 
deleted $d$ number of half lines. The number of connected components is 
equal to $1 + d$. By passing a multiple point of multiplicity $i$, then by 
Lemma \ref{topological-bifurcation}, we see that the number of connected components of 
$M_{\leq c}$ (resp. $M_{< c}$) increases exactly by $(i-1)$. Thus, after passing all multiple 
points, the number of connected components of $M_{< c}$, which is homeomorphic 
to $M({\mathcal A})$, is given by $1 + d + \sum_{i=2}^d (i-1)t_i$. 

\

\noindent
{\it Proof of Theorem \ref{general-dimensional-affine-diffeomorphism}.}
For a $c \in \R$ with $c \ll 0$, 
$M_{\leq c}$ (resp. $M_{< c}$) is diffeomorphic to the half space $\{ x_n \leq c \}$ 
(resp. $\{ x_n < c \}$ deleted $d$ number of half lines. 
By passing a multiple point of multiplicity $i$, for a sufficiently large $c$, 
$M_{\leq c}$ is obtained by 
attaching $i - 1$ number of trivial handles of index $n-2$, by Lemma \ref{topological-bifurcation}. 
After passing all multiple points, $M_{\leq c}$ is diffeomorphic to 
the space obtained by attaching $\sum_{i=2}^d (i-1)t_i$ number of trivial handles of index $n-2$ 
to the half space deleted $d$ number of half lines. 
Then $M_{< c}$ is diffeomorphic to the interior of $B_g$ with 
$g = d + \sum_{i=2}^d (i-1)t_i$. By Lemma  \ref{Topological-triviality}, 
for $c \in \R$ with $0 \ll c$, 
$M_{< c}$ is diffeomorphic to $M({\mathcal A})$. Hence we have Theorem \ref{general-dimensional-affine-topology}. 
\QED

\

\noindent
{\it Proofs of Theorem \ref{general-dimensional-affine-topology}
and Theorem \ref{three-dimensional-affine-topology}.} 
Theorem \ref{general-dimensional-affine-topology}
follows from Theorem \ref{general-dimensional-affine-diffeomorphism} and 
Theorem \ref{three-dimensional-affine-topology} follows from 
Theorem \ref{general-dimensional-affine-topology} by setting $n = 3$. 
\QED

\

\ber
{\rm 
Let $X$ be a subset of $\R^n$ which is a union of finite number of closed 
line segments and half lines. 
Then similarly to the proof of Theorem \ref{three-dimensional-affine-topology} using 
Lemma \ref{topological-bifurcation}, we see that, if there exists a height function 
$h : \R^n \to \R$ satisfying (i)(ii) such that $h\vert_X : X \to \R$ has no local maximum, then 
the complement $\R^n \setminus X$ is diffeomorphic to the interior of $n$-ball with trivially attached $g$-handles of index $n-2$, for some $g$. 
If $X \subset \R^n$ is compact, then any height function 
has a maximum, so non-trivial attachments may occur. 
}
\enr

\ber
\label{knot-complement}
{\rm 
The knot complements have much information than line arrangement complements. 
For example, 
it is known that, for knots $K, K' \subset S^3$, if $S^3 \setminus K$ and $S^3 \setminus K'$ 
are homeomorphic, then the pairs $(S^3, K)$ and $(S^3, K')$ are homeomorphic 
(\cite{GL}). Taking account of it, consider $(\R^3, X)$ 
for a line arrangement ${\mathcal A} = \{ \ell_1, \dots, \ell_d\}$ in $\R^3$ and $X := \bigcup_{i=1}^d \ell_i \subset \R^3$ and its one-point compactification $(S^3, \overline{X})$. Then the  
complement $S^3 \setminus X$ is homeomorphic to $M({\mathcal A})$ and to 
$B_g$, which depends only on the number $g = d + \sum_{i=1}^d (i-1)t_i$, 
while $g$ does not determine the topological type of the pair $(S^3, \overline{X})$ in general. 
}
\enr

\section{Projective line and linear plane arrangements}

Let $\widetilde{\mathcal A} = \{ \widetilde{\ell_1}, \dots, \widetilde{\ell_2}, \dots, \widetilde{\ell_d}\}$ be a  real projective line arrangement in the projective space $\R P^n$
and let 
${\mathcal B} = \{ L_1, L_2, \dots, L_d \}$ be the real linear plane arrangement 
in $\R^{n+2}$ corresponding to $\widetilde{\mathcal A}$. 
Then the complement $M({\mathcal B})$ of ${\mathcal B}$ 
is homeomorphic to the link complement $S^n \cap M({\mathcal B})$ times $\R_{>0}$, 
where $S^n$ is a sphere in $\R^{n+1}$ centred at the origin.  
Moreover $S^n \cap M({\mathcal B})$ 
is a double cover of $M(\widetilde{\mathcal A})$ for the 
corresponding projective line arrangement $\widetilde{\mathcal A}$ in $\R P^n$. 

\begin{figure}[ht]
\begin{center}
\includegraphics[width=3truecm, height=3truecm, clip, 
]{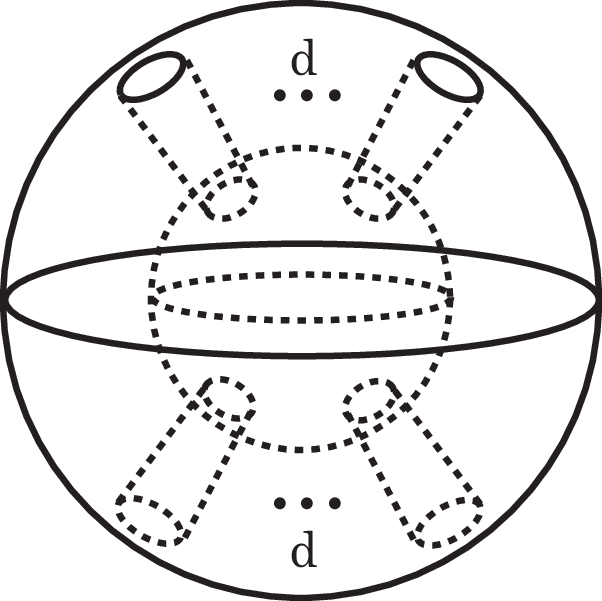} 
\caption{Punctured shell.}
\label{Punctured-shell}
\end{center}%
\end{figure}

Take a projective hyperplane $H \subset \R P^n$ such that 
$H$ intersects transversely to all lines $\widetilde{\ell_i}, 1 \leq i \leq d$, and that 
$H$ does not pass through any multiple point of $\widetilde{\mathcal A}$. 
Then identify $\R P^n \setminus H$ with the affine space $\R^n$ and the 
affine line arrangement ${\mathcal A}$ obtained by setting 
$\ell_i := \widetilde{\ell_i} \setminus H \subset \R^n$. 
Take a ball $D^n = \{ x \in \R^n \mid \Vert x\Vert \leq r \} \subset \R^n$ for a 
sufficiently large radius $r$ such that interior of $D^n$ contains all multiple points of 
${\mathcal A}$ and the boundary 
$\pa(D^n) = S^{n-1}$ intersects transversally to all lines $\ell_i, 1 \leq i \leq d$. 
Then the closure $\overline{U}$ of $U := \R P^n \setminus D^n$ is regarded as a tubular neighbourhood of $H$ in $\R P^n$. The closure $\overline{U}$ is homeomorphic to
the space $(S^{n-1} \times [-1, 1])/\!\sim$, where $(x, t) \sim (-x, -t)$. 
Let $a_1, \dots, a_{2d}$ be disjoint $2d$ points in $S^{n-1}$. Let 
$W_k^{n-1} \subset S^{n-1}$ be a sufficiently small 
open disk neighbourhood of $a_k, (1 \leq k \leq 2d)$. 
Set $N := S^{n-1} \setminus W_k^{n-1}$ and 
$\widetilde{N} :=  (N \times [-1, 1])/\!\sim \ (\subset \ (S^{n-1} \times [-1, 1])/\!\sim)$. 
Then $\widetilde{N}$ is an $n$-dimensional manifold with boundary $N$, 
which is double covered by a \lq\lq punctured shell" $N \times [-1, 1]$ 
(see Figure \ref{Punctured-shell}). 

Thus we observe 

\bep
\label{projective-arrangement} 
The intersection $U \cap M(\widetilde{\mathcal A})$ is homeomorphic to the interior of $\widetilde{N}$. 
The complement $M(\widetilde{\mathcal A}) \subset \R P^n$ is homeomorphic to the interior of 
$B_g \bigcup_\varphi \widetilde{N}$ for an attaching embedding 
$\varphi : N \to \pa(B_g)$. The homeomorphism class of $M(\widetilde{\mathcal A})$ is determined by the isotopy class of the embedding $\varphi$. The embedding $\varphi$ is determined by the intersection of $M({\mathcal A})$ and a hypersphere of sufficiently large radius in $\R^n$. 
\enp

\Proof
We see that 
the intersection of $M({\mathcal A})$ and a hypersphere of sufficiently large radius in $\R^n$ 
is homeomorphic to the sphere deleted $2d$-points. Then we have Proposition \ref{projective-arrangement} by Theorem \ref{general-dimensional-affine-topology}. 
\QED

\end{document}